\documentclass[11pt]{article}
\usepackage{amsmath,amssymb,amsbsy,amsfonts,amsthm,latexsym,
               amsopn,amstext,amsxtra,euscript,amscd,color}
\def \F {{\mathbb F}}

\def \Z {{\mathbb Z}}

\def \C {{\mathbb C}}

\def \Tr {{\rm Tr_n}}

\newtheorem{theorem}{Theorem}
\newtheorem{definition}[theorem]{Definition}
\newtheorem{lemma}[theorem]{Lemma}
\newtheorem{proposition}[theorem]{Proposition}
\newtheorem{remark}[theorem]{Remark}

\newtheorem{corollary}[theorem]{Corollary}

\def\aa{{\bf a}}

\def\cc{{\bf c}}
\def\uu{{\bf u}}

\def\xx{{\bf x}}
\def\yy{{\bf y}}
\def\zz{{\bf z}}

\def\00{{\bf 0}}
\def\11{{\bf 1}}
\def\+{\oplus}

\def \F {{\mathbb F}}

\def \Z {{\mathbb Z}}

\def \Tr {{\rm Tr_n}}

\begin{document}

\title{\huge\bf
\textrm{Modified planar functions and their components}}
\author{\Large  Nurdag\"ul Anbar$^{1}$, Wilfried Meidl$^{2}$,
\vspace{0.4cm} \\
\small $^1$Technical University of Denmark,\\ \small Matematiktorvet, Building 303B, DK-2800, Lyngby, Denmark\\ 
\small Email: {\tt nurdagulanbar2@gmail.com}\\
\small $^2$Johann Radon Institute for Computational and Applied Mathematics,\\
\small Austrian Academy of Sciences, Altenbergerstrasse 69, 4040-Linz, Austria\\
\small Email: {\tt meidlwilfried@gmail.com}
}

\date{}
\maketitle
\thispagestyle{empty}

\begin{abstract}
Zhou 2013 introduced modified planar functions to describe $(2^n,2^n,2^n,1)$ relative difference sets $R$ as a graph of a function on the finite 
field $\F_{2^n}$, and pointed out that projections of $R$ are difference sets that can be described by negabent or bent$_4$ functions, which are
Boolean functions given in multivariate form. Objective of this paper is to contribute to the understanding of these component functions of 
modified planar functions. We first completely describe a multivariate version of modified planar functions in terms of their bent$_4$ components.
In the second part we characterize the component functions of (univariate) modified planar functions in terms of appropriate generalizations of 
the Walsh-Hadamard transform, with respect to which they have a flat spectrum. We hereby obtain a description of modified planar functions by their
components which is similar to that of the classical planar functions in odd characteristic as a vectorial bent function.
\end{abstract}

\section{Introduction}

Let $G$ be a group of order $\mu\nu$ and let $N$ be a subgroup of $G$ of order $\nu$. A subset $R$ of $G$ of cardinality $k$ is called a 
$(\mu,\nu,k,\lambda)$-{\it relative difference set (RDS)} of $G$ relative to $N$, if every element of $G\setminus N$ can be written as a difference
of two elements of $R$ in exactly $\lambda$ ways, and there is no such representation for any nonzero element in $N$. The subgroup $N$
is hence also called the forbidden subgroup. 
A powerful description of RDSs is their description via characters (see for instance Section 2.4. in \cite{tpf}):
A subset $R$ of cardinality $k$ of a group $G$ of order $\mu\nu$ with  a subgroup $N$ of order $\nu$ is a $(\mu,\nu,k,\lambda)$-RDS of $G$ 
relative to $N$ if and only if
for every character $\chi$ of $G$
\begin{equation}
\label{RDSchar}
|\chi(R)|^2 = \left\{\begin{array}{ll}
                  k^2\quad & \mbox{if}\;\chi = \chi_0, \\
                  k-\lambda \nu\quad & \mbox{if}\;\chi \ne \chi_0,\:\mbox{and}\;\chi(g) = 1\;\mbox{for all}\;g\in N, \\
                  k\quad & \mbox{otherwise.}
                  \end{array} \right.
\end{equation}
As shown in \cite[Theorem 3.1]{gs}, an RDS relative to a normal subgroup $N$ of $G$ with parameters $(\mu,\nu,k,\lambda) = (q,q,q,1)$
uniquely describes a projective plane, hence is of particular interest. For abelian groups $G$ we have the following fundamental results
on the existence of $(q,q,q,1)$-RDSs:
\begin{itemize}
\item[I]: \cite{bjs} If an abelian group $G$ of odd order contains a $(q,q,q,1)$-RDS, then $q=p^n$ for some prime $p$ and $G$ contains an
elementary abelian subgroup of order $p^{n+1}$.
\item[II]: \cite{g} If an abelian group $G$ of even order contains a $(q,q,q,1)$-RDS, then $q=2^n$, $G\cong \Z_4^n$, and the forbidden subgroup is
$2\Z_4^n$.
\end{itemize}
{\bf I: $|G| = q^2$ and $q$ is odd.}

\noindent All known $(q,q,q,1)$-RDSs are subsets of the $2n$-dimensional vector space $G$ over the finite field $\F_p$. 
As a result we can represent the group $G$ as $G = (\F_p^n\times\F_p^n,+)$. In this case, $(q,q,q,1)$-RDSs can be expressed as the graph of a so-called {\it planar function} 
$f:\F_p^n\rightarrow\F_p^n$,
\[ R = \{(\xx,f(\xx))\,:\,\xx\in \F_p^n\}\ . \]
Planar functions $f$ are commonly defined by the property that the derivative in direction $\aa$
\[ D_{\aa}f(\xx) = f(\xx+\aa)-f(\xx) \]
is a permutation for all $\aa\ne \00$. Planar functions have been widely studied since their introduction in \cite{do}.

Observe that the group of characters of $(\F_p^n\times\F_p^n,+)$ is
\begin{align*}
\{\,\chi_{\uu,\cc}\,:\,\cc,\uu\in\F_p^n\,\} \quad\quad \text{with}\quad\quad \chi_{\uu,\cc}(\xx,\yy) := \epsilon_p^{\cc\cdot \yy-\uu\cdot \xx} \ ,
\end{align*}
where $\epsilon_p = e^{\frac{2\pi i}{p}}$ and $\xx\cdot\yy$ is the standard dot product of $\xx$ and $\yy$.
By the characterization of RDSs via characters, the set $R=\{(\xx,f(\xx))\,:\,\xx\in\F_p^n\}$ for $f:\F_p^n\rightarrow\F_p^n$
is an RDS of $(\F_p^n\times\F_p^n,+)$ (relative to the subgroup $\{\00\}\times\F_p^n$) if and only if
\begin{equation}
\label{pWalsh}
|\chi_{\uu,\cc}(R)|^2 = \left|\sum_{\xx\in\F_p^n}\epsilon_p^{\cc\cdot f(\xx)-\uu\cdot \xx}\right|^2 =: \left|\mathcal{W}_{\cc\cdot f}(\uu)\right|^2 = p^n
\end{equation}
for all $\uu \in \F_p^n$ and nonzero $\cc\in\F_p^n$. As is well known, the character sum in Equation $(\ref{pWalsh})$ is called the {\it Walsh-Hadamard transform} of the function 
$\cc\cdot f(\xx) =: f_{\cc}(\xx)$ from $\F_p^n$ to $\F_p$, which (for nonzero $\cc$) is called a {\it component function} of $f$. Moreover
$g:\F_p^n\rightarrow\F_p$ is called a {\it ($p$-ary) bent function} if $|\mathcal{W}_g(\uu)|^2 = p^n$ for every $\uu\in\F_p^n$, see \cite{ksw}.
We see that for all nonzero $\cc$, the component function $f_{\cc}$ of a planar function $f$ is bent. In fact we have the following
theorem, see e.g. \cite[Theorem 3.19]{p}.
\begin{theorem}
The function $f:\F_p^n\rightarrow\F_p^n$ is planar if and only if all its component functions $f_{\cc}$ are bent. 
\end{theorem}
For the equivalent representation of $G$ as $G = (\F_{p^n}\times\F_{p^n},+)$, the RDS can be represented as a graph of a univariate function 
$f:\F_{p^n}\rightarrow\F_{p^n}$. The Walsh transform of the component function $f_c(x) = \Tr(cf(x))$ is then given by
\[ \mathcal{W}_{f_c}(u) = \sum_{x\in\F_{p^n}}\epsilon_p^{\Tr(cf(x))-\Tr(ux)}, \]
where $\Tr(z)$ denotes the absolute trace of $z\in\F_{p^n}$.
All known planar functions are represented in univariate form and except from one example, they are all quadratic functions (Dembowski-Ostrom polynomials),
see \cite[Section 8]{p}.\\[.5em]
{\bf II: $|G| = q^2$ and $q$ is even.}

\noindent Recently, in the significant contribution \cite{z}, Zhou presented a solution how to describe a $(2^n,2^n,2^n,1)$-RDS in $\Z_4^n$ as a graph of a 
multivariate function $f$ from $\F_2^n$ to $\F_2^n$, respectively of a univariate function $f$ on $\F_{2^n}$.

In the multivariate case, the group $G \cong \Z_4^n$ is represented as $G = (\F_2^n\times\F_2^n, *)$ where
$(\xx_1,\yy_1) * (\xx_2,\yy_2) = (\xx_1+\xx_2,\yy_1+\yy_2+\xx_1\odot \xx_2)$ with 
$(x_1,x_2,\ldots,x_n) \odot (y_1,y_2,\ldots,y_n) = (x_1y_1,x_2y_2,\ldots,x_ny_n)$. 
The graph of a function $f:\F_2^n\rightarrow\F_2^n$ is then an RDS in $G$ relative to $\{\00\}\times\F_2^n \cong 2\Z_4^n$ if and only if
\begin{equation}
\label{modpla1}
f(\xx+\aa) + f(\xx) + \aa\odot \xx
\end{equation}
is a permutation for every nonzero $\aa\in\F_2^n$.

Alternatively we can represent $\Z_4^n$ as $G = (\F_{2^n}\times\F_{2^n}, \star)$ where the group operation is given by
$(x_1,y_1) \star (x_2,y_2) = (x_1+x_2,y_1+y_2+x_1x_2)$. In this case the graph of the univariate function $f:\F_{2^n}\rightarrow\F_{2^n}$ is an 
RDS in $G$ (relative to $\{0\}\times\F_{2^n}$) if and only if
\begin{equation}
\label{modpla2}
f(x+a) + f(x) + ax
\end{equation}
is a permutation for every nonzero $a\in\F_{2^n}$.

In accordance with \cite{p}, we call multivariate functions on $\F_2^n$ for which $(\ref{modpla1})$ is always a permutation, respectively univariate 
functions on $\F_{2^n}$ for which $(\ref{modpla2})$ is always a permutation, {\it modified planar functions}.

We emphasize that the multivariate modified planar functions for the group $ (\F_2^n\times\F_2^n, *)$ are not the one-to-one translation of the univariate
modified planar functions for the group $(\F_{2^n}\times\F_{2^n}, \star)$ by choosing an appropriate basis for $\F_{2^n}$. For details we refer to \cite{z},
where both versions were introduced. One can see that in detail they can behave quite differently, observing that every affine function from $\F_{2^n}$ to 
$\F_{2^n}$ is a trivial example of a univariate modified planar function. In particular, the zero-function on $\F_{2^n}$ is modified planar. This obviously 
is not the case for the zero function on $\F_2^n$. 

The component functions of modified planar functions correspond to $(2^n,2,2^n,2^{n-1})$-RDSs in $\Z_4\times\Z_2^{n-1}$, which essentially can be
represented by the graph of a {\it negabent function} (or more general a {\it bent$_4$} function, which is also called a {\it shifted bent} function), 
see the discussions in \cite[Section 5]{z} and \cite[Section 7]{psz}. One of the objectives of this paper is to contribute to the understanding of the 
component functions of modified planar functions.

Many interesting results on modified planar functions such as relations to semifields, which are analogues of the results on planar functions in odd 
characteristic, have been discovered. For more detailed information we refer to \cite{psz,sz,z}, and the excellent recent overview paper \cite{p}.

Objective of this paper is to contribute to the understanding of the component functions of modified planar functions. In Section \ref{multi},
we first recall bent$_4$ functions, which have been defined and analysed as multivariate Boolean functions with a flat spectrum with respect
to certain unitary transforms. Then we elaborate connections between (multivariate) modified planar functions and bent$_4$ functions, which were
already pointed at in \cite{z}, in more detail.
In Section \ref{uni} we characterize the component functions of univariate modified planar functions in terms of character sums, again unitary 
transforms, with respect to which they have a flat spectrum. Though they behave somewhat different than the multivariate bent$_4$ functions, we
suggest to call them (univariate) bent$_4$ functions as well.

\section{Components of multivariate modified planar functions}
\label{multi}

For the convenience of the reader we first set up some notation. We denote by $\C$ the set of complex numbers, by $\F_{2^n}$ the finite field of order $2^n$ and 
by $\F_2^n$ the space of all $n$-tuples $(x_1,\ldots,x_n)$ of elements from $\F_2$. In particular, we 
denote by $\mathbf{0}$ and $\mathbf{1}$ the $n$-tuples in $\F_{2}^n$ whose entries are all $0$ and $1$, respectively.

For $\cc = (c_1,\ldots,c_n)$, $\xx = (x_1,\ldots,x_n)$ in $\F_2^n$, in accordance with the notation in \cite{gps} we set
\[ s_2^\cc(\xx) := \bigoplus_{1\le i<j\le n}(c_ix_i)(c_jx_j) \ . \]
Note that $s_2^\cc(\xx) = s_2(\cc\odot \xx)$, where $s_2(x)$
is the homogeneous symmetric Boolean function with algebraic degree $2$. For an element $\cc\in\F_2^n$ and a Boolean function $g:\F_2^n\rightarrow\F_2$, a unitary 
transform $\mathcal{U}_g^\cc:\F_2^n\rightarrow\C$ is defined by (cf.\cite{gps})
\begin{equation}
\label{UTF1}
\mathcal{U}_g^\cc(\uu) = \sum_{\xx\in\F_2^n}(-1)^{g(\xx)+ s_2^\cc(\xx)}i^{\cc\cdot \xx}(-1)^{\uu\cdot \xx}.
\end{equation}
By \cite[Proposition 3]{gps}, for $\cc,\xx\in\F_2^n$ we have
\[ \cc\cdot \xx + 2s_2^\cc(\xx) \equiv wt(\cc\odot \xx) \bmod 4 \]
where $wt(\zz)$ denotes the Hamming weight of $\zz\in\F_2^n$. As a result we can write Equation $(\ref{UTF1})$ also as
\begin{equation}
\label{UTF2}
\mathcal{U}_g^\cc(\uu) = \sum_{\xx\in\F_2^n}(-1)^{g(\xx)+ \uu\cdot \xx}i^{wt(\cc\odot \xx)} \ .
\end{equation}
Obviously, for $\cc=\mathbf{0}$ the transform $\mathcal{U}_g^{\mathbf{0}}$ in Equation $(\ref{UTF2})$ is the conventional Walsh-Hadamard transform of $g$, and for $\cc = \mathbf{1}$ it  
is known as the {\it nega-Hadamard transform} of $g$, see \cite{pp}.
\begin{definition}
\label{bent4.1}
A Boolean function $g:\F_2^n\rightarrow\F_2$ is called bent$_4$ if it has a flat spectrum with respect to at least
one of the transforms $\mathcal{U}_g^\cc$. In other words, $g$ is bent$_4$ if there exists an element $\cc\in\F_2^n$ such that
$|\mathcal{U}_g^\cc(u)| = 2^{n/2}$ for all $\uu\in\F_2^n$.
\end{definition}
The Boolean functions with flat spectrum with respect to $\mathcal{U}_g^{\mathbf{0}}$ are the celebrated bent functions.
If $g$ is flat with respect to $\mathcal{U}_g^{\mathbf{1}}$, then $g$ is called {\it negabent}. 
We will call $g$ a $\cc$-bent$_4$ function if $g$ has a flat spectrum with respect to $\mathcal{U}_g^\cc$.
Negabent functions have been investigated in \cite{pp,spp,sgcgm,spt,zwp}, for more background on bent$_4$ functions we refer to \cite{rp}, 
where they have been introduced, and to \cite{gps}.

In this section we look at the relation between modified planar functions in multivariate form defined as in $(\ref{modpla1})$
and bent$_4$ functions in more detail. 
\begin{proposition}
Let $G = (\F_2^n\times\F_2^n, *)$ with $$(\xx_1,\yy_1) * (\xx_2,\yy_2) = (\xx_1+\xx_2,\yy_1+\yy_2+\xx_1\odot \xx_2) \ . $$ 
The group of characters of $G$ is then $\chi_{G} = \{\chi_{\uu,\cc}\,:\,\uu,\cc\in\F_2^n\}$ where
\[ \chi_{\uu,\cc}(\xx,\yy) = (-1)^{\uu\cdot \xx + \cc\cdot \yy}i^{wt(\cc\odot \xx)} \ . \]
\end{proposition}  
{\it Proof.}
We first show that $\chi_{\uu,\cc}:G\rightarrow\C$ is a group homomorphism. For $\uu,\cc\in\F_2^n$ we have
\begin{align*}
\chi_{\uu,\cc}((\xx_1,\yy_1)*(\xx_2,\yy_2)) & = \chi_{\uu,\cc}(\xx_1+\xx_2, \yy_1+\yy_2+\xx_1\odot \xx_2) \\
&= (-1)^{\uu\cdot(\xx_1+\xx_2)+\cc\cdot(\yy_1+\yy_2+\xx_1\odot \xx_2)}i^{wt(\cc\odot (\xx_1+\xx_2))} \ . 
\end{align*}
On the other hand,
\[ \chi_{\uu,\cc}(\xx_1,\yy_1)\chi_{\uu,\cc}(\xx_2,\yy_2) = (-1)^{\uu\cdot(\xx_1+\xx_2)+\cc\cdot(\yy_1+\yy_2)}i^{wt(\cc\odot \xx_1)+wt(\cc\odot \xx_2)}. \]
Hence we conclude that $\chi_{\uu,\cc}((\xx_1,\yy_1)* (\xx_2,\yy_2)) = \chi_{\uu,\cc}(\xx_1,\yy_1)\chi_{\uu,\cc}(\xx_2,\yy_2)$ if and only if
\[ i^{2\cc\cdot(\xx_1\odot \xx_2)+wt(\cc\odot (\xx_1+\xx_2))} = i^{wt(\cc\odot \xx_1) + wt(\cc\odot \xx_2)} \ , \]
or equivalently
\begin{equation}
\label{mod4} 
2\cc\cdot(\xx_1\odot \xx_2)+wt(\cc\odot (\xx_1+\xx_2)) \equiv wt(\cc\odot \xx_1) + wt(\cc\odot \xx_2) \bmod 4. 
\end{equation}
With 
\begin{align*}
wt(\cc\odot(\xx_1+\xx_2))& = wt(\cc\odot \xx_1 + \cc\odot \xx_2) \\
&= wt(\cc\odot \xx_1) + wt(\cc\odot\xx_2) - 2wt(\cc\odot \xx_1\odot \xx_2)\ ,
\end{align*}
we see that Equation $(\ref{mod4})$ holds if and only if
\begin{align*}
2\cc\cdot(\xx_1\odot\xx_2) \equiv 2wt(\cc\odot \xx_1\odot \xx_2) \bmod 4 \ ;
\end{align*}
that is $\cc\cdot(\xx_1\odot\xx_2) \equiv wt(\cc\odot \xx_1\odot \xx_2) \bmod 2$,
which trivially holds.

It remains to show that $\chi(\uu_1,\cc_1) \ne \chi(\uu_2,\cc_2)$ if $(\uu_1,\cc_1) \ne (\uu_2,\cc_2)$. Suppose that $\chi(\uu_1,\cc_1)=\chi(\uu_2,\cc_2)$,
which yields
\begin{equation}
\label{C=C}
(-1)^{(\uu_1+\uu_2)\cdot \xx + (\cc_1+\cc_2)\cdot \yy}i^{wt(\cc_1\odot \xx)-wt(\cc_2\odot \xx)} = 1 
\end{equation}
for all $\xx,\yy\in\F_2^n$. If $(\ref{C=C})$ holds, then $wt(\cc_1\odot \xx)-wt(\cc_2\odot \xx) \equiv 0\bmod 2$.
Let $\cc_1 = (c_{1,1},\ldots,c_{1,n})$, $\cc_2 = (c_{2,1},\ldots,c_{2,n})$, suppose that w.l.o.g., $c_{1,1} \ne c_{2,1}$, 
and let $\xx = (1,0,\ldots,0)$. Then $wt(\cc_1\odot \xx)-wt(\cc_2\odot \xx)=\pm 1$.
We conclude that Equation $(\ref{C=C})$ implies $\cc_1=\cc_2$. Equation $(\ref{C=C})$ then implies $(-1)^{(\uu_1+\uu_2)\cdot \xx}=1$ for all $\xx\in\F_2^n$,
thus $\uu_1 = \uu_2$. \hfill$\Box$\\[.5em]
Let $f$ be a function on $\F_2^n$. By the characterization of RDSs via characters in Equation $(\ref{RDSchar})$, the set
$\{(\xx,f(\xx))\,:\,\xx\in\F_2^n\}$ is a $(2^n,2^n,2^n,1)$-RDS in $G$ if and only if for every $\cc,\uu\in\F_2^n$, $\cc\ne \00$,
\[ \left|\sum_{\xx\in\F_2^n}\chi_{\uu,\cc}(\xx,f(\xx))\right| = \left|\sum_{\xx\in\F_2^n}(-1)^{\cc\cdot f(\xx)+ \uu\cdot \xx}i^{wt(\cc \odot \xx)}\right| =
|\mathcal{U}_{f_\cc}^\cc(\uu)| = 2^{n/2}, \] 
where $f_\cc$ denotes the component function $f_\cc(\xx) = \cc\cdot f(\xx)$ of $f$.
Note that $|\sum_{\xx\in\F_2^n}\chi_{\00,\00}(\xx,f(\xx))| = 2^{n}$ and $|\sum_{\xx\in\F_2^n}\chi_{\uu,\00}(\xx,f(\xx))| = 0$
if $\uu\ne \00$, trivially hold. We obtain the following theorem.
\begin{theorem}
\label{mpmu}
The function $f:\F_2^n\rightarrow\F_2^n$ is a modified planar function if and only if for every nonzero
$\cc\in\F_2^n$ the component function $f_\cc$ is a $\cc$-bent$_4$ function.
\end{theorem}
{\it Proof.} By the above discussion, both conditions are equivalent to $\{(\xx,f(\xx))\,:\,\xx\in\F_2^n\}$ being
a $(2^n,2^n,2^n,1)$-RDS in $G$. \hfill$\Box$\\[.5em]
The equivalence of the conditions 
\begin{itemize}
\item[(i)] $f(\xx)+f(\xx+\aa)+\aa\odot \xx$ is a permutation for every nonzero $\aa$, and 
\item[(ii)] for every nonzero $\cc$ the component function $f_\cc$ is $\cc$-bent$_4$,
\end{itemize}
can of course also be deduced directly. As for the analog conditions for planar functions in odd characteristic, by squaring the respective transform
one first shows that a function $g:\F_2^n\rightarrow\F_2$ is $\cc$-bent$_4$ if and only if
\begin{equation}
\label{bal4} 
g(\xx)+g(\xx+\zz) + \cc\cdot(\zz\odot\xx)
\end{equation}
is balanced for all nonzero $\zz\in\F_2^n$. We omit this argument at this position and include it in the next section, where we deal with modified 
planar functions in univariate form.

\section{Components of univariate modified planar functions}
\label{uni}

So far, all known examples of modified planar functions are in univariate representation; i.e. they are univariate
functions $f$ satisfying Equation $(\ref{modpla2})$, hence induce RDSs in the group $G = (\F_{2^n}\times\F_{2^n}, \star)$ 
with the operation $(x_1,y_1) \star (x_2,y_2) = (x_1+x_2,y_1+y_2+x_1x_2)$.
Trivial examples are all affine functions. Besides from  affine functions, all known examples of modified planar functions are 
quadratic functions (represented by Dembowski-Ostrom polynomials). The existence of a nonquadratic modified planar function is 
an open problem. 

To analyse component functions of modified planar functions on $\F_{2^n}$, we first determine the group of characters 
of $G$. For $c,x\in\F_{2^n}$ we define
\[ \sigma(c,x) := \sum_{0\le i<j\le n-1}(cx)^{2^i}(cx)^{2^j}. \]
Note that $\sigma(c,x)^2 = \sigma(c,x)$, and hence $\sigma(c,x)$ is a Boolean function.
\begin{lemma}
\label{s-prpty}
For $c,x_1,x_2\in\F_{2^n}$ we have
\[ \sigma(c,x_1+x_2) = \sigma(c,x_1) + \sigma(c,x_2) + \Tr(cx_1)\Tr(cx_2) + \Tr(c^2x_1x_2). \]
\end{lemma}
{\it Proof.}
Expanding $\sigma(c,x_1+x_2)$, we get the desired conclusion as follows.
\begin{align*}
&\sigma(c,x_1+x_2)\\
 &= \sum_{i<j}(c(x_1+x_2))^{2^i}(c(x_1+x_2))^{2^j} \\
&= \sum_{i<j}(cx_1)^{2^i}(cx_1)^{2^j} + \sum_{i<j}(cx_2)^{2^i}(cx_2)^{2^j} + \sum_{i<j}((cx_1)^{2^i}(cx_2)^{2^j} + (cx_1)^{2^j}(cx_2)^{2^i}) \\
&= \sigma(c,x_1) + \sigma(c,x_2) + \sum_{i,j}(cx_1)^{2^i}(cx_2)^{2^j} + \sum_{i}(cx_1)^{2^i}(cx_2)^{2^i} \\
& = \sigma(c,x_1)+  \sigma(c,x_2) + \left(\sum_{i}(cx_1)^{2^i}\right)\left(\sum_{j}(cx_2)^{2^j}\right) + \sum_{i}(c^2x_1x_2)^{2^i} \\
&= \sigma(c,x_1) +\sigma(c,x_2) + \Tr(cx_1)\Tr(cx_2) + \Tr(c^2x_1x_2).
\end{align*}
\hfill$\Box$ 
\begin{proposition}
Let $G = (\F_{2^n}\times\F_{2^n},\star)$ be the group with the operation $(x_1,y_1) \star (x_2,y_2) = (x_1+x_2,y_1+y_2+x_1x_2)$.
Then the group of characters of $G$ is $\chi_{G} = \{\chi_{u,c}\,:\,u,c\in\F_{2^n}\}$, where
\[ \chi_{u,c}(x,y) = (-1)^{\Tr(ux)+ \Tr(c^2y) + \sigma(c,x)}i^{\Tr(cx)} \ . \]
\end{proposition}
{\it Proof.}
With the definition of $\chi_{u,c}$, we see that the equality 
\[ \chi_{u,c}((x_1,y_1)\star (x_2,y_2)) = \chi_{u,c}(x_1,y_1)\chi_{u,c}(x_2,y_2) \] 
holds if and only if
\[ (-1)^{\sigma(c,x_1) + \sigma(c,x_2)}i^{\Tr(cx_1)+\Tr(cx_2)} = (-1)^{\Tr(c^2x_1x_2)+\sigma(c,x_1+x_2)}i^{\Tr(c(x_1+x_2))} \ ,\]
or equivalently
\begin{align}
\label{kjfj}
\nonumber
& 2\sigma(c,x_1)+2\sigma(c,x_2)+\Tr(cx_1)+\Tr(cx_2) \equiv 2\Tr(c^2x_1x_2)+2\sigma(c,x_1+x_2) \\
& +\Tr(c(x_1+x_2)) \bmod 4. 
\end{align}
With the identity
\begin{equation}
\label{identity} 
\Tr(x) + \Tr(y) \equiv \Tr(x+y)+2\Tr(x)\Tr(y) \bmod 4, 
\end{equation}
Condition $(\ref{kjfj})$ is equivalent to
\[ 2\sigma(c,x_1)+2\sigma(c,x_2)+2\Tr(cx_1)\Tr(cx_2) \equiv 2\Tr(c^2x_1x_2)+2\sigma(c,x_1+x_2) \bmod 4 \ ; \]
i.e. $\sigma(c,x_1) + \sigma(c,x_2) + \Tr(cx_1)\Tr(cx_2) + \Tr(c^2x_1x_2) \equiv \sigma(c,x_1+x_2) \bmod 2$,
which is true by Lemma \ref{s-prpty}.

It remains to show that $\chi_{u_1,c_1} \ne \chi_{u_2,c_2}$ if $(u_1,c_1) \ne (u_2,c_2)$. 
Suppose that $\chi_{u_1,c_1} = \chi_{u_2,c_2}$; i.e.
\begin{equation}
\label{chi=chi}
(-1)^{\Tr((u_1-u_2)x)+ \Tr((c_1^2-c_2^2)y)+ \sigma(c_1x)+ \sigma(c_2x)}i^{\Tr(c_1x)-\Tr(c_2x)} = 1
\end{equation}
for all $x\in \F_{2^n}$. Observe that (modulo $4$), $\Tr(c_1x)-\Tr(c_2x) \in \{-1,0,1\}$. As a result, Equation $(\ref{chi=chi})$ implies 
that $\Tr(c_1x)-\Tr(c_2x) = 0$ for all $x\in\F_{2^n}$, which implies $c_1=c_2$. Immediately one sees that then $u_1=u_2$. 

\hfill$\Box$\\[.5em]
Knowing the group of characters, we can employ the character sum characterization of RDSs in Equation $(\ref{RDSchar})$ to obtain conditions for functions
$f:\F_{2^n}\rightarrow\F_{2^n}$ for which the graph of $f$ is an RDS in $G$.
Observing that $|\sum_{x\in\F_{2^n}}\chi_{0,0}(x,f(x))| = 2^{n}$ and $|\sum_{x\in\F_{2^n}}\chi_{u,0}(x,f(x))| = 0$
if $u\ne 0$, are again trivially satisfied, the set $\{(x,f(x))\,:\,x\in\F_{2^n}\}$ is an RDS in $G$ if and only if for every 
$u,c\in\F_{2^n}$, $c\ne 0$, we have
\[ \left|\sum_{x\in\F_{2^n}}\chi_{u,c}(x,f(x))\right| = \left|\sum_{x\in\F_{2^n}}(-1)^{\Tr(c^2f(x))+ \sigma(c,x)}i^{\Tr(cx)}(-1)^{\Tr(ux)}\right| = 2^{n/2}. \]
This motivates the definition of a unitary transform $\mathcal{V}_{g}^c$, $c\in\F_{2^n}$, for a function 
$g:\F_{2^n}\rightarrow\F_2$ as
\[ \mathcal{V}_{g}^c(u) := \sum_{x\in\F_{2^n}}(-1)^{g(x)+ \sigma(c,x)}i^{\Tr(cx)}(-1)^{\Tr(ux)} \ . \]
Note that for $c=0$, we again obtain the conventional Walsh-Hadamard transform for univariate Boolean functions.
We define the nega-Hadamard transform for univariate functions as $\mathcal{V}_{g}^c$ for $c=1$.

As the graph of $f$ is a $(2^n,2^n,2^n,1)$-RDS in $G$ if and only if for all nonzero $c\in\F_{2^n}$
the component function $f_{c^2}(x) = \Tr(c^2f(x))$ of $f$ has a flat spectrum with respect to $\mathcal{V}_g^c$, we 
find it natural to define bent$_4$ functions in univariate form by using these transforms.
\begin{definition}
\label{bent4.2}
A function $g:\F_{2^n}\rightarrow\F_2$ is called bent$_4$ if it has a flat spectrum with respect to at least
one of the transforms $\mathcal{V}_g^c$. In other words, $g$ is bent$_4$ if there exists an element $c\in\F_{2^n}$ such that
$|\mathcal{V}_g^c(u)| = 2^{n/2}$ for all $u\in\F_{2^n}$.
\end{definition}
The functions from $\F_{2^n}$ to $\F_2$ with a flat spectrum with respect to $\mathcal{V}_g^0$ are the bent functions. 
In accordance with the case of multivariate functions we call a function with a flat spectrum with respect to 
$\mathcal{V}_g^1$ a negabent function.
With this notations we obtain a univariate analog of Theorem \ref{mpmu}.
\begin{theorem}
\label{mpun}
A function $f:\F_{2^n}\rightarrow\F_{2^n}$ is a modified planar function if and only if for every nonzero
$c\in\F_{2^n}$ the component function $f_{c^2}(x)=\Tr(c^2f(x))$ is a $c$-bent$_4$ function.
\end{theorem}
Note that $c\rightarrow c^2$ is a permutation of the multiplicative group of $\F_{2^n}$, and hence Theorem \ref{mpun} gives a condition on 
all component functions of $f$.

We finally show the equivalence of the condition in Theorem \ref{mpun} and Condition $(\ref{modpla2})$
directly, which may also provide us a further understanding of component functions of modified planar functions.
As for the case of conventional bent functions we follow the approach via Hadamard matrices, see e.g. \cite{d}.
\begin{lemma}
\label{Fcon}
Let $h$ be a complex valued function on $\F_{2^n}$. Then
\[ \Psi(u) = \sum_{z\in\F_{2^n}}h(z)(-1)^{\Tr(uz)+\sigma(c,z)}i^{\Tr(cz)} = h(0) \]
for all $u\in \F_{2^n}$, if and only if $h(z) = 0$ for $z\ne 0$.
\end{lemma}
\begin{theorem}
A function $g:\F_{2^n}\rightarrow\F_2$ is a $c$-bent$_4$ function if and only if
\[ g(x) + g(x+z) +\Tr(c^2xz) \]
is balanced for all $z\ne 0$.
\end{theorem}
{\it Proof.} 
For every $u\in\F_{2^n}$ we have
\begin{equation}
\label{V^2}
\mathcal{V}_g^c(u)\overline{\mathcal{V}_g^c(u)} = \sum_{z\in\F_{2^n}}(-1)^{\Tr(uz)+\sigma(c,z)}i^{\Tr(cz)}
\sum_{x\in\F_{2^n}}(-1)^{g(x)+ g(x+z)+\Tr(c^2xz)} \ ,
\end{equation}
where $\overline{\mathcal{V}_g^c(u)}$ is the complex conjugate of $\mathcal{V}_g^c(u)$. Hence we conclude that if $g(x)+g(x+z)+\Tr(c^2xz)$ is balanced, 
then $\mathcal{V}_g^c(u)\overline{\mathcal{V}_g^c(u)} = 2^n$; i.e. $g$ is $c$-bent$_4$.

Conversely suppose that $g$ is $c$-bent$_4$, thus $\mathcal{V}_g^c(u)\overline{\mathcal{V}_g^c(u)} = 2^n$.
Setting
\[ h(z) := \sum_{x\in\F_{2^n}}(-1)^{g(x) + g(x+z) +\Tr(c^2xz)} \ , \]
with Equation $(\ref{V^2})$ this yields
\[ \sum_{z\in\F_{2^n}}h(z)(-1)^{\Tr(uz)+ \sigma(c,z)}i^{\Tr(cz)} = h(0) = 2^n \ . \]
By Lemma \ref{Fcon}, this implies that $h(z) = 0$ for every nonzero $z$, which holds if and only if
$g(x) + g(x+z) + \Tr(c^2xz)$ is balanced for every $z\ne 0$.
\hfill$\Box$\\[.5em]
The function $h:\F_{2^n}\rightarrow\F_{2^n}$ is a permutation if and only if for every nonzero $c \in \F_2^n$, $\Tr(ch)$ is balanced (see \cite[Theorem 7.7]{ln}).
This yields the desired result given in the following corollary.
\begin{corollary}
Let $f$ be a function on $\F_{2^n}$. Then $\Tr(c^2f)$ is $c$-bent$_4$ for all nonzero $c\in\F_{2^n}$ if and only if
$f(x)+f(x+z)+xz$ is a permutation for all nonzero $z\in\F_{2^n}$
\end{corollary}
\begin{remark}
\label{unimul}
Univariate bent$_4$ functions defined as in Definition \ref{bent4.2} behave somewhat different than the multivariate 
bent$_4$ functions in Definition \ref{bent4.1}. Observe that every affine function $f(x) = \Tr(ax)+b$, $a\in\F_{2^n}$,
$b\in\F_2$, is $c$-bent$_4$ for every nonzero $c$. Affine functions in multivariate form are trivial examples for negabent 
functions (see also \cite[Proposition 1]{pp}). But a multivariate affine function is never $\cc$-bent$_4$ for any $\cc\ne\11$:
As easily seen, Equation $(\ref{bal4})$ is constant for an affine function $g$ if one chooses $\zz=(z_1,\ldots,z_n)$ such that $z_i=0$
if $c_i=1$.
\end{remark}
\begin{remark}
In \cite{sar,zq} a negabent function $g:\F_{2^n}\rightarrow\F_2$ is defined to be a function for which
$g(x) + g(x+z) + \Tr(xz)$ is balanced for all nonzero $z$. This is equivalent with our definition of 
a negabent function as a function with a flat spectrum with respect to $\mathcal{V}_g^1$.
\end{remark}

\section{Acknowledgment}
Nurdag\"{u}l Anbar gratefully acknowledges the support from The Danish Council for Independent Research (Grant No. DFF--4002-00367) and H.C. \O rsted COFUND Post-doc Fellowship from the project ``Algebraic curves with many rational points".\\ 
Wilfried Meidl is supported by the Austrian Science Fund (FWF) Project no. M 1767-N26.


\begin{thebibliography}{99}


\bibitem{bjs} A. Blokhuis, D. Jungnickel, B. Schmidt, Proof of the prime power conjecture for projective planes of order $n$ with abelian collineation groups of order $n^2$. 
Proc. Amer. Math. Soc. 130 (2002), no. 5, 1473--1476.




\bibitem{do} P. Dembowski, T.G. Ostrom, Planes of order $n$ with collineation groups of order $n^2$. Math. Z. 103 (1968) 239--258. 

\bibitem{d} J.F. Dillon, Elementary Hadamard difference sets, Ph.D. dissertation, University of Maryland, 1974.
%
\bibitem{gps} S. Gangopadhyay, E. Pasalic, P. St\u anic\u a, A note on generalized bent criteria for Boolean functions. 
IEEE Trans. Inform. Theory 59 (2013), no. 5, 3233--3236. 

\bibitem{gs} M.J. Ganley, E. Spence, Relative difference sets and quasiregular collineation groups. J. Combinatorial Theory Ser. A 19 (1975), no. 2, 134--153.

\bibitem{g} M.J. Ganley, On a paper of P. Dembowski and T. G. Ostrom: "Planes of order $n$ with collineation groups of order $n^2$''. Arch. Math. (Basel) 27 (1976), no. 1, 93--98.

%

\bibitem{ksw} P.V. Kumar, R.A. Scholtz, L.R. Welch, Generalized bent functions and their properties. J. Combin. Theory Ser. A 40 (1985), no. 1, 90--107.

\bibitem{ln} R. Lidl, H. Niederreiter, Finite Fields, 2nd ed., Encyclopedia Math. Appl., vol. 20, Cambridge Univ. Press, Cambridge, (1997).
%

\bibitem{pp} M.G. Parker, A. Pott, On Boolean functions which are bent and negabent. Sequences, subsequences, and consequences, 9--23, 
Lecture Notes in Comput. Sci., 4893, Springer, Berlin, 2007. 

\bibitem{psz} A. Pott, K.U. Schmidt, Y. Zhou, Semifields, relative difference sets, and bent functions. Algebraic curves and finite fields, 161--178, 
Radon Ser. Comput. Appl. Math., 16, De Gruyter, Berlin, 2014.

\bibitem{p} A. Pott, Almost perfect and planar functions, Des. Codes Cryptogr. 78 (2016), 141--195.

\bibitem{rp} C. Riera, M.G. Parker, Generalized bent criteria for Boolean functions. I. IEEE Trans. Inform. Theory 52 (2006), no. 9, 4142--4159.


\bibitem{sar} S. Sarkar, Characterizing negabent Boolean functions over finite fields. Sequences and their applications—SETA 2012, 77–88, 
Lecture Notes in Comput. Sci., 7280, Springer, Heidelberg, 2012. 

\bibitem{spp} K.U. Schmidt, M.G. Parker, A. Pott, Negabent functions in the Maiorana-McFarland class. Sequences and their applications—SETA 2008, 390–402, 
Lecture Notes in Comput. Sci., 5203, Springer, Berlin, 2008.

\bibitem{sz} K.U. Schmidt, Y. Zhou, Planar functions over fields of characteristic two. J. Algebraic Combin. 40 (2014), no. 2, 503--526.

\bibitem{sgcgm} P. St\u anic\u a, S. Gangopadhyay, A. Chaturvedi, A.K. Gangopadhyay, S. Maitra, Investigations on bent and negabent functions via the 
nega-Hadamard transform, IEEE Trans. Inform. Theory 58 (2012), 4064--4072.

\bibitem{spt} W. Su, A. Pott, X. Tang, Characterization of negabent functions and construction of bent-negabent functions with maximum algebraic degree. 
IEEE Trans. Inform. Theory 59 (2013), 3387--3395.  

\bibitem{tpf} Y. Tan, A. Pott, T. Feng, Strongly regular graphs associated with ternary bent functions. J. Combin. Theory Ser. A 117 (2010), no. 6, 668--682.

\bibitem{zwp} F. Zhang, Y. Wei, E. Pasalic, Constructions of bent-negabent functions and their relation to the completed Maiorana-McFarland class. 
IEEE Trans. Inform. Theory 61 (2015), 1496--1506.

\bibitem{z} Y. Zhou, $(2n,2n,2n,1)$-relative difference sets and their representations. J. Combin. Des. 21 (2013), no. 12, 563--584.

\bibitem{zq} Y. Zhou, L. Qu, Constructions of negabent functions over finite fields, Cryprography and Communications, to appear.

\end{thebibliography}
\end{document}